\documentclass[12pt]{article}
\usepackage[margin=25mm,a4paper]{geometry}
\usepackage{amsmath,amssymb,amsthm}
\usepackage[vcentermath]{youngtab}

\begin{document}

\newcommand{\GL}{{\mathrm{GL}}}
\newcommand{\SL}{{\mathrm{SL}}}
\newcommand{\Gr}{{\mathrm{Gr}}}
\newcommand{\Hom}{{\mathrm{Hom}}}
\newcommand{\codim}{{\mathrm{codim}}}
\newcommand{\cF}{{\mathcal{F}}}
\newcommand{\cB}{{\mathcal{B}}}
\newcommand{\cO}{{\mathcal{O}}}
\newcommand{\ssn}{\subsetneq}
\newcommand{\PP}{{\mathbb{P}}}
\newcommand{\NN}{{\mathbb{N}}}
\newcommand{\ZZ}{{\mathbb{Z}}}
\newcommand{\CC}{{\mathbb{C}}}
\renewcommand{\AA}{{\mathbb{A}}}
\newcommand{\KK}{{\mathbb{K}}}
\newcommand{\Supp}{{\mathrm{Supp}\;}}
\newcommand{\Span}{{\mathrm{Span}\;}}
\newcommand{\Fl}{{\mathrm{Fl}\;}}
\newcommand{\Id}{{\mathrm{Id}}}
\newcommand{\diag}{{\mathrm{diag}\;}}
\newcommand{\rk}{{\mathrm{rk}\;}}
\newcommand{\ybullet}{$\bullet$}

\newtheorem{theorem}{Theorem}
\newtheorem{prop}[theorem]{Proposition}
\newtheorem{corollary}[theorem]{Corollary}
\newtheorem{lemma}[theorem]{Lemma}

\newcommand{\defin}{{\noindent\bf Definition. }}
\newcommand{\remark}{{\noindent\bf Remark. }}
\newcommand{\example}{{\noindent\bf Example. }}

\title{Desingularizations of Schubert varieties in double Grassmannians}
\author{Evgeny Smirnov}
\date{}
\maketitle

\abstract{Let $X=\Gr(k,V)\times\Gr(l,V)$ be the direct product of
two Grassmann varieties of $k$- and $l$-planes in a
finite-dimensional vector space $V$, and let $B\subset \GL(V)$ be
the isotropy group of a complete flag in $V$. We consider $B$-orbits
in $X$, which are an analog to Schubert cells in Grassmannians. We
describe this set of orbits combinatorially and construct
desingularizations for the closures of these orbits, similar to the
Bott--Samelson desingularizations for Schubert varieties. }

\section{Introduction}

Let $V$ be a finite-dimensional vector space over an arbitrary ground field $\KK$.

We are interested in describing pairs of subspaces in $V$ of fixed
dimensions $k$ and $l$ up to the action of the group
$B\subset\GL(V)$ of non-degenerate upper-triangular matrices. So,
what we describe is the set of $B$-orbits in the direct product
$X=\Gr(k,V)\times \Gr(l,V)$ of two Grassmann varieties. This
decomposition is similar to the Schubert decomposition for
Grassmannians or to the Ehresmann--Bruhat decomposition for full
flag varieties.

The combinatorial description of $B$-orbits in $X$ was given (as a
special case of a more general problem) by Magyar, Weyman and
Zelevinsky \cite{MWZ}; their proof substantially uses quiver theory.
The description given below does not refer to these results and is
only based some elementary linear algebra. This is a generalization
of the description, given in Pin's thesis \cite{Pin}, of orbits in
the symmetric space $\GL_{k+l}/(\GL_k\times \GL_l)$.

We are also interested in the closures of $B$-orbits in $X$. They
can be viewed as analogues of Schubert varieties in Grassmannians.
The singularities of Schubert varieties are well-known objects: they
admit nice desingularizations, constructed by Bott and Samelson; it
is known that Schubert varieties are normal and have rational
singularities; their singular loci can be described explicitly.
(Good references on this topic are the lecture notes \cite{B2} and
the textbook \cite{Man}.) So it is natural to ask the same questions
(desingularization, normality, and rationality of singularities) for
the closures of $B$-orbits in $X$. In this paper we construct
desingularizations of these varieties.

Our interest in this problem is also motivated by the recent paper
\cite{BZ} by G.~Bobi\'nski and G.~Zwara, where they prove that the
singularities of the closures of orbits in representations of
 type $D$ quivers are equivalent to the singularities of Schubert
varieties in double Grassmannians.

The author is grateful to Michel Brion for constant attention to
this work, and to Ernest Vinberg and Dmitri Timashev for useful
remarks and comments.

\section{Description of orbits}

\subsection{Notation}

Let $V$ be an $n$-dimensional vector space over a field $\KK$, and
let  $k,l<n$ be positive integers. The results of Sec.~2 are valid
over an arbitrary ground field; however, in Secs.~3 and 4 we assume
$\KK$ to be algebraically closed. The direct product
$\Gr(k,V)\times\Gr(l,V)$ is denoted by $X$. Usually we do not
distinguish between points of $X$ and the corresponding
configurations $(U,W)$ of subspaces, where $U,W\subset V$, $\dim
U=k$, and $\dim W=l$.

We fix a Borel subgroup $B$ in $\GL(V)$. Let
$V_\bullet=(V_1,\dots,V_n=V)$ be the complete flag in $V$ stabilized
by $B$.

\subsection{Combinatorial description}

In this section, we introduce some combinatorial objects that
parametrize pairs of subspaces up to $B$-action. Namely, pairs of
subspaces are parametrized by triples consisting of two Young
diagrams contained in the rectangles of size $k\times (n-k)$ and
$l\times(n-l)$, respectively, and an involutive permutation in
$S_n$.

Along with constructing these combinatorial objects, we also
construct some ``canonical'' bases in the spaces $U$, $W$, and $V$.

\begin{prop}\label{baseschoice}
(i). There exist ordered bases $(u_1,\dots,u_k)$, $(w_1,\dots,w_l)$,
and $(v_1,\dots,v_n)$ of $U$, $W$, and $V$, respectively, such that:

\begin{itemize}

\item $V_i=\langle v_1,\dots,v_i\rangle$ for each $i\in\{1,\dots,n\}$ (angle brackets stand for the linear span of
vectors).

\item $u_i=v_{\alpha_i}$, where $i\in\{1,\dots,k\}$ and $\{\alpha_1,\dots, \alpha_k\}\subset
\{1,\dots,n\}$.

\item The $w_i$ are either basis vectors of $V$ or vectors with
two-element ``support'': $w_i=v_{\beta_i}$ or
$w_i=v_{\gamma_i}+v_{\delta_i}$, where $\gamma_i>\delta_i$;
moreover, in the latter case $v_{\gamma_i}\in U$ (that is,
$\{\gamma_1,\dots,\gamma_r\}\subset\{\alpha_1,\dots,\alpha_k\}$).

\item All the $\beta_i$, $\gamma_i$ and $\delta_i$ are distinct;
moreover, every $\delta_i$ is distinct from every $\alpha_j$.

\end{itemize}

(ii). With the notation of (i), define a permutation $\sigma\in S_n$
as the product of all the transpositions $(\gamma_i,\delta_i)$.
Their supports are disjoint, and so this product does not depend of
their order.

 Then, for a given pair $(U,W)$, the subsets
$\bar\alpha=\{\alpha_1,\dots,\alpha_k\}$,
$\bar\beta=\{\beta_1,\dots,\beta_{l-r}\}$,
$\bar\gamma=\{\gamma_1,\dots,\gamma_{r}\}$ of $\{1,\dots, n\}$, and
the permutation $\sigma$ are independent of the choice of bases in
$U$, $W$, and $V$.
\end{prop}

\begin{proof} (i) We prove this by induction over $n$.

If $n=1$, there is nothing to prove.

The inductive step from $n-1$ to $n$ goes as follows. Take a nonzero
vector $v_1\in V_1$, and consider the following cases:

\begin{itemize}

\item $v_1\notin U+W$. Let $\bar V=V/\langle v_1\rangle$. Let us take the flag $\bar V_\bullet=(\bar
V_2\subset\dots\subset \bar V_n)$ in this quotient space and
consider the images $(\bar U,\bar W)$ of $(U,W)$ under taking the
quotient. Both these subspaces are isomorphic to their preimages,
$\bar U\cong U$ and $\bar W\cong W$. Now we apply the induction
hypothesis to this configuration as follows. Let us choose ordered
bases $\{\bar u_1,\dots,\bar u_k\}$, $\{\bar w_1,\dots,\bar w_l\}$,
and $\{\bar v_1,\dots,\bar v_{n-1}\}$ in $\bar U$, $\bar W$, and
$\bar V$. Then we choose a lift  $\imath\colon\bar V\hookrightarrow
V$. Now take the preimages of these basis vectors in $V$ in the
following way: $u_i=\imath(\bar u_i)$, $w_i=\imath(\bar w_i)$, and
$v_i=\imath(\bar v_{i-1})$. We obtain the desired triple of bases.

\item $v_1\in U$, $v_1\notin W$. Set $u_1=v_1$ and again
apply the induction hypothesis to the quotient $\bar V=V/\langle
v_1\rangle$ with the flag $\bar V_\bullet$ and the configuration
$(\bar U,\bar W)$. The only difference is that in this case $\dim
\bar U=\dim U-1$. After that we lift the bases from
$\bar U$, $\bar W$, and $\bar V$ to $V$ in a similar way.

\item The case when $v_1\notin U$ and $v_1\in W$ is similar to the
previous one. (We set $w_1=v_1$.)

\item If $v_1\in U\cap W$, we set $u_1=w_1=v_1$ and again apply
the induction hypothesis.

\item The most interesting case is the last one: $v_1\in U+W$, but $v_1\notin U$ and $v_1\notin W$. Consider then the set
of vectors $S=\{v\mid v\in U, v_1+v\in W\}$. Since $v_1$ belongs to
the sum $U+W$, it follows that this set is nonempty. Now let $j$ be
the minimum number such that $V_j$ contains vectors in $S$, and that
$v_j\in V_j\cap S$. Set $u_1=v_j$ and $w_1=v_1+v_j$. Now apply the
induction hypothesis to the $(n-2)$-dimensional space $\bar
V=V/\langle v_1,v_j\rangle$, the configuration of two subspaces
$\bar U=U/\langle v_1,v_j\rangle$, $\bar W=W/\langle
v_1,v_j\rangle$, and the flag
\begin{multline*} \bar V_\bullet =(
V_2/V_1\subset\dots\subset V_{j-1}/V_1=\\=V_j/\langle v_1,v_j\rangle
\subset V_{j+1}/\langle v_1, v_j\rangle \subset\dots\subset
V_n/\langle v_1,v_j\rangle).
\end{multline*}
We lift basis vectors from $\bar V$ to $V$ as follows:
$$
v_i=\imath(\bar v_{i-1}),\text{ if }i\in[2,j-1];\qquad
v_i=\imath(\bar v_{i-2})\text{ if }i\in[j+1,n],
$$
where, as above, $\imath$ is an embedding of $\bar V$ in $V$. We
have already defined the vectors $v_1$ and $v_j$.
\end{itemize}

(ii) Take a configuration $(U,W)$ and assume that there exist two
triples of ordered bases, $((u_1,\dots,u_k),$ $(w_1,\dots,w_l),$
$(v_1,\dots,v_n))$ and $((u'_1,\dots,u'_k),$ $(w'_1,\dots,w'_l),$
$(v'_1,\dots,v'_n))$, satisfying the conditions of (i) and such that
either the triples of sets $(\bar\alpha,\bar\beta, \bar\gamma)$ and
$(\bar\alpha',\bar\beta', \bar\gamma')$ or the permutations $\sigma$
and $\sigma'$ corresponding to the first and the second triple of
bases, respectively, are not equal.

The set $\bar\alpha$ can be described as follows: $i\in\bar\alpha$
if and only if $\dim U\cap V_i >\dim U\cap V_{i-1}$. This means that
$\bar\alpha=\bar\alpha'$.

By the same argument, we can prove that $\bar\beta\cup\bar\gamma
=\bar\beta'\cup\bar\gamma'$.

Now let us prove that $\sigma=\sigma'$. This will complete the
proof, since
$\bar\beta=\{j\in\bar\beta\cup\bar\gamma\mid\sigma(j)=j\}$.

Let $j$ be the minimum number in $\bar\beta\cup\bar\gamma$ such that
$\sigma(j)\neq\sigma'(j)$. Suppose that $\sigma(j)<\sigma'(j)$. Two
cases may occur:

a) $i:=\sigma'(j)\neq j$. First, observe that $i\notin\bar\alpha$.
Consider the subspace
\begin{multline*} \tilde V=(U\cap V_j)+V_{i-1}=\langle
v_s,v_{\alpha_i}\mid s\leq
i-1,\alpha_i\in\bar\alpha\cup[i,j]\rangle=\\ = \langle
v'_s,v'_{\alpha_i}\mid s\leq
i-1,\alpha_i\in\bar\alpha'\cup[i,j]\rangle.
\end{multline*}

Let $R=\{r\in\bar\beta\cup\bar\gamma\mid
r,\sigma(r)\in[1,i-1]\cup(\bar\alpha\cap[i,j])\}$ and $R'
=\{r\in\bar\beta\cup\bar\gamma\mid
r,\sigma'(r)\in[1,i-1]\cup(\bar\alpha\cap[i,j])\}$. One can readily
see that
$$
\dim \tilde V\cap W=\#R=\#R'.
$$
But $\sigma(r)=\sigma'(r)$ for all $r\in[1,j-1]$, and $j$ belongs to
$R$ and does not belong to $R'$. That means that the cardinalities
of these two sets are different, that gives us the desired
contradiction.

b) If $\sigma'(j)=j$, set $i=\sigma(j)$, and proceed as in a).
\end{proof}

Now introduce a combinatorial construction that parametrizes
configuration types. Namely, given a configuration, we will
construct a pair of Young diagrams with some boxes marked.

Suppose we have a configuration $(U,W)$. Take the sets $\bar\alpha$,
$\bar\beta$, $\bar\gamma$, and the involution $\sigma$ corresponding
to this configuration, as described in Prop.~\ref{baseschoice}.
Consider a rectangle of size $k\times(n-k)$ and construct a path
from its bottom-left to top-right corner such that each $j$th step
is vertical if $j$ belongs to $\bar\alpha$ (that is, $v_j$ is equal
to some $u_i$) and horizontal otherwise. This path bounds (from
below) the first Young diagram.

The second diagram is contained in a rectangle of size
$l\times(n-l)$. Again, we construct a path bounding it. Let the
$j$th step of this path be vertical if $j\in\bar\beta\cup\bar
\gamma$ and horizontal otherwise.

If  $j\in\bar\gamma$, then the $\sigma(j)$th step of this path is
horizontal. This also means that the $j$th and $\sigma(j)$th steps
of the path bounding the first diagram are vertical and horizontal,
respectively. In each diagram, take the box located above the
$\sigma(j)$th step and to the left of the $j$th step and put a dot
into each of these boxes.

We refer to this pair of diagrams with dots as a \emph{marked pair}.

\example Let $n=9$, $k=4$, and $l=3$. Suppose that
$\bar\alpha=\{3,5,6,9\}$, $\bar\beta=\{2,5\}$, $\bar\gamma=\{9\}$,
$\sigma=(7,9)$. Then the corresponding marked pair of diagrams is as
follows:
$$
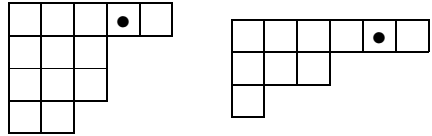
\young(\hfil\hfil\hfil\bullet\hfil,\hfil\hfil\hfil,\hfil\hfil\hfil,\hfil\hfil)
\qquad\young(\hfil\hfil\hfil\hfil\bullet\hfil,\hfil\hfil\hfil,\hfil)
$$

\remark Note that the constructed diagrams (without dots) are the
same as the diagrams that correspond to the Schubert cells
containing the points $U\in\Gr(k,V)$ and $W\in\Gr(l,V)$. (The
correspondence between Schubert cells and Young diagrams is
described, for example, in the textbooks \cite{F} and \cite{Man}.)

\subsection{Stabilizers and dimensions of orbits}

Now let us find the stabilizer $B_{(U,W)}$ for a given configuration
$(U,W)$.

\begin{prop}\label{isotr} In the notation of Prop.~ref{baseschoice}, the
stabilizer of a configuration $(U,W)$ written in the basis
$(v_1,\dots,v_n)$ consists of the upper-triangular matrices
$A=(a_{ij})\in\GL(n)$ satisfying the following conditions:
\begin{enumerate}

\item[(a)] $a_{\gamma\gamma}=a_{\sigma(\gamma)\sigma(\gamma)}$ for each
$\gamma\in\bar\gamma$.

\item[(b)] $a_{i \alpha}=0$ for all $\alpha\in \bar\alpha$ and $i\notin
\bar\alpha$.

\item[(c)] $a_{j\beta}=0$ for all $\beta\in \bar\beta$ and
$j\notin \bar\beta\cup \bar\gamma\cup\sigma(\bar\gamma)$.

\item[(d)]  $a_{\gamma\beta}=a_{\sigma(\gamma)\beta}$ for all $\beta\in \bar\beta$ and $\gamma\in \bar\gamma$,
$\gamma<\beta$.

\item[(e)] $a_{j\gamma}=-a_{j\sigma(\gamma)}$ for all $j\notin \bar\beta\cup\bar\gamma\cup\sigma(\bar\gamma)$ and $\gamma\in
\bar\gamma$.

\item[(f)] for any $\gamma_1,\gamma_2\in \bar\gamma$,
$\gamma_1<\gamma_2$, one of the following cases occurs:

\begin{itemize}

\item $\sigma(\gamma_2)<\sigma(\gamma_1)<\gamma_1<\gamma_2$; then
$a_{\gamma_1\gamma_2}=a_{\sigma(\gamma_1)\gamma_2}=a_{\sigma(\gamma_2)\gamma_1}=a_{\sigma(\gamma_1)\sigma(\gamma_2)}=0$.

\item $\sigma(\gamma_1)<\sigma(\gamma_2)<\gamma_1<\gamma_2$; then
$a_{\sigma(\gamma_2)\gamma_1}=a_{\sigma(\gamma_1)\gamma_2}=0$,
$a_{\gamma_1\gamma_2}=a_{\sigma(\gamma_1)\sigma(\gamma_2)}$.

\item $\sigma(\gamma_1)<\gamma_1<\sigma(\gamma_2)<\gamma_2$; then
$a_{\sigma(\gamma_1)\gamma_2}=0$,
$a_{\gamma_1\gamma_2}+a_{\gamma_1\sigma(\gamma_2)}=a_{\sigma(\gamma_1)
\sigma(\gamma_2)}$.

\end{itemize}

\end{enumerate}
\end{prop}

\begin{corollary}
The stabilizer of a configuration $(U,W)$ is a semidirect product of
a toric and a unipotent part,
$$
B_{(U,W)}=T_{(U,W)}\ltimes U_{(U,W)},
$$
where $T_{(U,W)}$ is the subgroup determined by conditions (a) in
the group of nondegenerate diagonal matrices, so that $\dim
T_{(U,W)}=n-\#\bar\gamma$, and $U_{(U,W)}$ is the subgroup
determined by conditions (b)--(f) in the group of unitriangular
matrices.
\end{corollary}

\defin The codimension of the toric part of the stabilizer is
called the \emph{rank} of a configuration (or its corresponding
$B$-orbit):
$$
\rk (U,W):=n-\dim T_{(U,W)}=\#\bar\gamma.
$$

\begin{proof}[Proof of Proposition~\ref{isotr}] First, the stabilizer
$B_{(U,W)}$ is a subgroup of $B$ and hence  consists of
upper-triangular matrices.

Next, it preserves the subspace $U=\langle
v_{\alpha_1},\dots,v_{\alpha_k}\rangle$. This means that each
transformation $A\in B_{(U,W)}$ takes each $v_{\alpha_i}$ to a
linear combination of $v_{\alpha_j}$, so that all the elements $a_{i
\alpha }$ whith $\alpha\in \bar\alpha$ and $i\notin \bar\alpha$ are
zero. (Note that the zeros in $A$ themselves form the Young diagram
corresponding to the subspace $U$, rotated $90^\circ$ clockwise.
This proves, in particular, that the dimension of a Schubert cell in
a Grassmannian is equal to the number of boxes in the corresponding
Young diagram.)

Thus, the boxes of the first Young diagram are in a one-to-one
correspondence with the linear equations defining $B_{U}$ as a
subgroup of the group of upper-triangular matrices: the box located
above the $i$th (horizontal) step and to the left of the $j$th
(vertical) step  of the corresponding path (denote this box by
$(i,j)$) corresponds to the equation $a_{ij}=0$.

Likewise, the stabilizer of our configuration preserves the subspace
$W$. This gives another set of linear equations on the entries
$a_{ij}$, and the number of these equations is equal to the number
of boxes in the second diagram of the corresponding marked pair.
Again, we establish a one-to-one correspondence between the boxes of
this diagram and these equations, denoting boxes as in the previous
paragraph. This correspondence is as follows:
\begin{itemize}

\item $a_{j\beta}=0$ for all $\beta\in \bar\beta$ and
$j\notin \bar\beta\cup \bar\gamma\cup\sigma(\bar\gamma)$, $j<\beta$.
The corresponding box is $(j,\beta)$.

\item $a_{j\gamma}=-a_{j\sigma(\gamma)}$ for all $j\notin \bar\beta\cup\bar\gamma\cup\sigma(\bar\gamma)$
and $\gamma\in \bar\gamma$, $j<\gamma$. The corresponding box is
$(j,\gamma)$.

\item $a_{\sigma(\gamma)\gamma}+a_{\gamma\gamma}-a_{\sigma(\gamma)\sigma(\gamma)}=0$
for all $\gamma\in\bar\gamma$. The corresponding box is
$(\sigma(\gamma),\gamma)$.

\item $a_{\gamma\beta}=a_{\sigma(\gamma)\beta}$ for all $\beta\in \bar\beta$ and $\gamma\in \bar\gamma$,
$\gamma<\beta$. The corresponding box is $(\sigma(\gamma),\beta)$.

\item
$a_{\sigma(\gamma_1)\sigma(\gamma_2)}+a_{\sigma(\gamma_1)\gamma_2}=a_{\gamma_1\sigma(\gamma_2)}+a_{\gamma_1\gamma_2}$
for any $\gamma_1<\gamma_2$. This equation corresponds to the box
$(\sigma(\gamma_1),\gamma_2)$.

\end{itemize}

Bringing all these equations together completes the proof of the
proposition.
\end{proof}

Once we know the stabilizer of a configuration, we can calculate its
dimension and hence the dimension of the orbit $B(U,W)\subset X$.
Analyzing the equations above, one can find a combinatorial
interpretation of dimension in terms of marked pairs of Young
diagrams.

To do this, we have to introduce one more combinatorial notion.
Suppose we have two rectangles of size $k\times(n-k)$ and $l\times
(n-l)$, respectively, and a path in each of these rectangles
bounding a Young diagram. (Both paths are of length $n$.) Consider
the set of all numbers $i$ such that the $i$th steps in the paths
bounding both diagrams are horizontal and take the columns in the
diagrams lying above these steps. After that do the same for the
pairs of steps that are ``simultaneously vertical'' and take the
rows to the left of these steps.

The intersection of columns and rows we have taken is a Young
diagram as well. We refer to it as a \emph{common diagram}
corresponding to the given pair of diagrams.

\example The pair $(Y_1,Y_2)$ of Young diagrams
$$
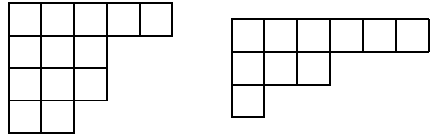
\yng(5,3,3,2)\qquad \yng(6,3,1)
$$
has the following common diagram $Y_{com}$:
$$
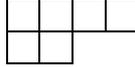
\yng(4,2)
$$

By our construction of marked pairs, dots can only be contained in
boxes of the common diagram of a marked pair.

\begin{corollary} Let $(U,W)$ be a configuration of subspaces,
and let $(Y_1,Y_2)$ be the corresponding marked pair of Young
diagrams, with dots in some boxes of the common diagram $Y_{com}$.

Consider the diagram $Y_{com}$. For a given box, take the set formed
by all the boxes in the same row above it, all the boxes in the same
row to the left of it, and this box itself. Such set is said to be
an \emph{inward hook} with vertex in the initial box. Now take all
its boxes with dots in $Y_{com}$ and consider all inward hooks with
vertices in these boxes. Let $H$ be the set of boxes that belong to
at least one of these inward hooks.  Then the dimension of the
$B$-orbit of $(U,W)$ is equal to
$$
\dim B(U,W)=\#Y_1+\#Y_2-\#Y_{com}+\#H,
$$
where $\#Y$ is the number of boxes in $Y$.
\end{corollary}

\remark $\#H$ is equal to the total number of boxes contained in all
hooks, rather than the sum lengths of all hooks. This means that a
box included into two hooks should be counted once, not twice.

\begin{proof} In the proof of Prop.~\ref{isotr}, we deal with two systems
of linear equations on the matrix entries $a_{ij}$. They correspond
to the stabilizers of the subspaces $U$ and $W$ and consist of
$\#Y_1$ and $\#Y_2$ equations, respectively. One can readily see
that the equations corresponding to the box $(i,j)$ coincide in both
systems if the box $(i,j)$ of the common diagram does not belong to
any hook, and also that the system of equations obtained as the
union of these two systems is linearly independent. Thus, the
codimension of $B_{(U,W)}$ in $B$ (that is, the dimension of
$B(U,W)$) is equal to $\#Y_1+\#Y_2-\#Y_{com}+\#H$.
\end{proof}

\example Let the common diagram for a marked pair be as follows:
$$
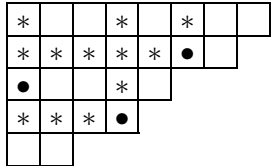
\young(*\hfil\hfil*\hfil*\hfil\hfil,*****\bullet\hfil,\bullet\hfil
\hfil*\hfil,***\bullet,\hfil\hfil)
$$
Then $\#Y_{com}=26$ and $\#H=15$ (Here $H$ consists of all
nonempty boxes.)

In particular, the dimension formula allows us to describe the
minimal (maximally degenerate) and maximal (open) orbit.  The
minimal orbit is zero-dimensional and corresponds to
$Y_1=Y_2=\varnothing$. It is the point $(\langle
v_1,\dots,v_k\rangle, \langle v_1,\dots,v_l\rangle)\in X$. Both
Young diagrams corresponding to the maximal orbit are rectangular,
of size $k\times(n-k)$ and $l\times(n-l)$, respectively. Their
common diagram is also a rectangle of size
$\min\{k,l\}\times(n-\max\{k,l\})$, with dots on the diagonal
issuing from the bottom-right corner.

{\sc Example.} For $n=8$, $k=3$, and $l=4$, the combinatorial data
corresponding to the maximal orbit are as follows:

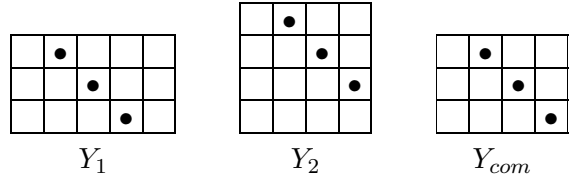
\begin{center}
\begin{tabular}{ccccc}
\young(\hfil\bullet\hfil\hfil\hfil,\hfil\hfil\bullet\hfil\hfil,\hfil\hfil\hfil\bullet\hfil)&\qquad&
\young(\hfil\bullet\hfil\hfil,\hfil\hfil\bullet\hfil,\hfil\hfil\hfil\bullet,\hfil\hfil\hfil\hfil)&\qquad&
\young(\hfil\bullet\hfil\hfil,\hfil\hfil\bullet\hfil,\hfil\hfil\hfil\bullet)
\\
$Y_1$ && $Y_2$&& $Y_{com}$
\end{tabular}
\end{center}

\subsection{Decomposition of $X$ into a union of $\GL(V)$-orbits}

The $\GL(V)$-orbits in $X$ have a much simpler description; they are
parametrized by one nonnegative integer, namely, by the dimension of
the intersection $U\cap W$. For this number (denote it by $i$) we
have the inequality
$$
\max\{0,k+l-n\}\leq i\leq \min\{k,l\}.
$$
Denote the corresponding $\GL(V)$-orbit by $X_i$:
$$
X=\bigsqcup\limits_{i\in\{\max(0,k+l-n),\dots,\min(k,l)\}} X_i.
$$
It follows from our construction of the combinatorial data
corresponding to a $B$-orbit that $\dim(U\cap W)=\#(\bar\alpha\cap
\bar\beta)$.

\section{Weak order on the set of orbits}

In the previous section, we have described the set of $B$-orbits in
$\Gr(k,V)\times\Gr(l,V)$. There exist several partial order
structures on this set. The first, most natural one is defined as
follows:

\defin Let $\cO$ and $\cO'$ be two $B$-orbits in
$\Gr(k,V)\times\Gr(l,V)$. We say that $\cO$ is less than or equal to
$\cO'$ with respect to the \emph{strong} (or \emph{topological})
\emph{order} if and only if $\cO\subset\bar\cO'$. (Here and further,
a bar stands for the Zariski closure of a set). Notation:
$\cO\le\cO'$.

There exists another order on this set, usually called the weak
order. Here the notation  is taken from \cite{B1}.

Let $W$ be the Weyl group for $\GL(n)$, and let $\Delta$ be the
corresponding root system. Denote the simple reflections by
$s_1,\dots,s_{n-1}$, and the corresponding simple roots by
$\alpha_1,\dots,\alpha_{n-1}$. Let $P_i=B\cup B s_i B$ be the
minimal parabolic subgroup in $\GL(V)$ corresponding to the simple
root $\alpha_i$.

We say that $\alpha_i$ \emph{raises} an orbit $\cO$ to $\cO'$, if
$\bar\cO'=P_i\bar\cO\neq\bar\cO$. In this case,
$\dim\cO'=\dim\cO+1$. This notion allows us to define the weak
order.

\defin An orbit $\cO$ is said to be less than or equal to  $\cO'$
with respect to the \emph{weak order} (notation: $\cO\preceq\cO'$)
if $\bar\cO'$ can be obtained as the result of several consecutive
raisings of $\bar\cO$ by minimal parabolic subgroups:
$$
\cO\preceq\cO'\quad\Leftrightarrow\quad \exists
(i_1,\dots,i_r)\colon \bar\cO'=P_{i_r}\dots P_{i_1}\bar\cO.
$$

Let us represent this relation of order by an oriented graph.
Consider a graph $\Gamma(X)$ with vertices indexed by $B$-orbits in
$X$. Join $\cO$ and $\cO'$ with an edge labeled by $i$ and directed
to $\cO'$, if $P_i$ raises $\cO$ to $\cO'$.

It is clear that each connected component of $\Gamma(X)$ consists of
$B$-orbits contained in a same $\GL(V)$-orbit $X_i$ and  has a
unique maximal element, namely, the $B$-orbit that is open in $X_i$.

Our next aim is to describe minimal elements with respect to the
weak order in each connected component.

\subsection{Combinatorial description of minimal parabolic subgroup
action}\label{combdescr}

Consider an orbit $\cO$ and the corresponding combinatorial data:
the sets $\bar\alpha$, $\bar\beta$ and $\bar\gamma$ and the
involution $\sigma\in S_n$. Let a minimal parabolic subgroup
$P_i=B\cup Bs_i B$ raise the orbit $\cO$ to an orbit $\cO'\neq\cO$.
Now we shall describe the combinatorial data ($\bar\alpha'$,
$\bar\beta'$, $\bar\gamma'$, $\sigma'$) of $\cO'$.

Denote the transposition $(i,i+1)\in S_n$ by $\tau_i$.

The following cases may occur:

\begin{enumerate}

\item
Suppose that
$$
i\in \bar\alpha,\quad i\notin\bar\beta, \quad i+1\notin \bar\alpha,
\quad i+1\in\bar\beta,
$$
or, on the opposite,
$$
i\notin\bar\alpha,\quad i\in\bar\beta, \quad i+1\in \bar\alpha,\quad
i+1\notin\bar\beta.
$$

These two cases correspond to two orbits that can be raised to
$\cO'$ by $P_i$. In this case, the new combinatorial data for $\cO'$
is given as follows:
$$
\begin{array}{rcl}
\bar\alpha' & =& \bar\alpha\cup\{i+1\}\setminus\{i\};\\
\bar\beta' & =& \bar\beta\setminus\{i,i+1\};\\
\bar\gamma' & =& \bar\gamma\cup\{i+1\}\\
\sigma'&=&\sigma\cdot\tau_i.
\end{array}
$$

Note that $\rk\cO'=\rk\cO+1$ and $\dim\cO'=\dim\cO+1$.

In the language of marked pairs of diagrams, this is represented as
follows. If the $i$th and the $i+1$th steps of the path bounding the
first diagram form an indentation (i.e., the $i$th step is vertical
and the $(i+1)$st one is horizontal) and the corresponding intervals
of the second diagram form a spike (or, on the opposite, we have a
spike in the first diagram and an indentation in the second), then
both these pairs of steps are replaced by spikes bounding a marked
box.

\example Apply the minimal parabolic subgroup $P_2$ to the orbit
$\cO\subset\Gr(3,7)\times\Gr(4,7)$ defined by the  marked pair
$$
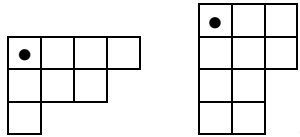
\text{\young(\bullet\hfil\hfil\hfil,\hfil\hfil\hfil,\hfil)}\qquad
\text{\young(\bullet\hfil\hfil,\hfil\hfil\hfil,\hfil\hfil,\hfil\hfil)}.
$$
The orbit $\cO'$ obtained as the result of this raising is defined
by the marked pair
$$
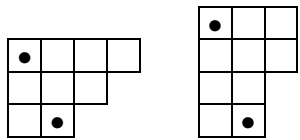
\text{\young(\bullet\hfil\hfil\hfil,\hfil\hfil\hfil,\hfil\bullet)}\qquad
\text{\young(\bullet\hfil\hfil,\hfil\hfil\hfil,\hfil\hfil,\hfil\bullet)}.
$$

\item In all other cases, $\bar\alpha'=
\tau_i(\bar\alpha)$, $\bar\beta'=\tau_i(\bar\beta)$, and
$\bar\gamma'=\tau_i(\bar\gamma)$, and the permutation $\tilde\sigma$
is the result of the conjugation of $\sigma$ by $\tau_i$:
$$
\tilde\sigma=\tau_i\sigma\tau_i.
$$
The ranks of these orbits are equal, $\rk\cO'=\rk\cO$.

\end{enumerate}

\subsection{Minimal orbits}

\begin{lemma}
All minimal $B$-orbits with respect to the weak order in a given
$\GL(V)$-orbit have rank 0.
\end{lemma}

\proof Assume the converse. Let $\cO$ be a minimal orbit with
nonzero rank, and let $(\bar\alpha,\bar\beta,\bar\gamma,\sigma)$ be
the corresponding combinatorial data such that $\sigma\ne Id$. Let
$p\in \bar\gamma$ and  $p'=\sigma(p)$. Without loss of generality we
can suppose that there is no other $q\in \bar\gamma$ such that
$p<q<\sigma(q')<p'$.

Let $C_1$ denote the set of spikes in the first diagram between
$p$th and $p'$th steps, that is, the set of indices $i$, $p\leq
i<p'$, such that the $i$th step in the first diagram is horizontal
and the $(i+1)$st is vertical. Similarly, let $D_1$ denote the set
of indentations, that is, the set of $i$, $p\leq i<p'$, such that
the $i$th step is vertical and the $(i+1)$st is horizontal.

Denote the similar sets for the second diagram by $C_2$ and $D_2$.
Note that $\#C_1=\#D_1+1$ and $\#C_2=\#D_2+1$, since in both
diagrams the $p$th steps are horizontal and the $p'$th steps are
vertical.

Now take a $j\in(C_1\setminus D_2)\cup(C_2\setminus D_1)$. Using
argument in \ref{combdescr}, we can show that there exists an orbit
$\cO'$ such that $\bar\cO=P_j\bar\cO'$. We describe the
combinatorial data for this orbit.

If the permutation $\sigma$ contains the transposition $(j,j+1)$,
then the combinatorial data for $\cO'$ is as follows:
$$
\begin{array}{rcl}
\bar\alpha' & =& \bar\alpha\cup\{j\}\setminus\{j+1\};\\
\bar\beta' & =& \bar\beta\cup\{j\};\\
\bar\gamma' & =& \bar\gamma\setminus\{j+1\}\\
\sigma'&=&\sigma\cdot\tau_j.
\end{array}
$$
Otherwise, $\bar\alpha'= \tau_j(\bar\alpha)$,
$\bar\beta'=\tau_j(\bar\beta)$, $\bar\gamma'=\tau_j(\bar\gamma)$,
and $\sigma'=\tau_j\sigma\tau_j$.

The calculation of the dimensions shows that $\dim\cO'=\dim\cO-1$.

To complete the proof, it remains to show that the set
$(C_1\setminus D_2)\cup(C_2\setminus D_1)$ is nonempty:
\begin{multline*}
\#((C_1\setminus D_2)\cup(C_2\setminus D_1))\geq\max
(\#(C_1\setminus D_2),
 \#(C_2\setminus D_1))\geq\\
\geq \max(\#C_1-\#C_2+1,\#C_2-\#C_1+1)\geq 1.
\end{multline*}
\endproof

After that we can find all the minimal orbits in $X_d$. One can
readily see that each minimal orbit has the following combinatorial
data:
$$
\begin{array}{rcl}
\bar\alpha\cup\bar\beta & =& \{1,\dots,k+l-d\};\\
\bar\alpha\cap\bar\beta & =& \{1,\dots,d\};\\
\bar\gamma & =& \varnothing;\\
\sigma & =& Id.\\
\end{array}
$$
The dimensions of all minimal orbits in $X_d$ are equal to
$(k-d)(l-d)$. In particular, this means that they all are closed in
$X_d$, since $X_d$ does not contain any orbits of lower dimension.
They correspond to the decompositions of the set
$\{d+1,\dots,k+l-d\}$ into two parts, $\bar\alpha\setminus\bar\beta$
and $\bar\beta\setminus\bar\alpha$, and so their number is equal to
$\binom{k+l-2d}{k-d}$.

Also note that the pair of Young diagrams that corresponds to a
minimal orbit is complementary: one can put these two diagrams
together so that they fill a rectangle of size
$(k-d)\times(l-d)$.

The common diagram corresponding to each of these pairs of diagrams
is empty. This means that all minimal orbits are stable under
the $(B\times B)$-action, that is, they are direct products of two
Schubert cells in two Grassmannians.

These results can be summarized as the following theorem.

\begin{theorem}\label{minorbits}
Each $X_d$, where $d\in\{\max(k+l-n,0),\dots,\min(k,l)\}$, contains
$\binom{k+l-2d}{k-d}$ minimal orbits. All these orbits are closed in
$X_d$ and have dimension $(k-d)(l-d)$. They are direct products of
Schubert cells.
\end{theorem}

\section{Desingularizations of the orbit closures}

In this section, we construct desingularizations for the closures of
$B$-orbits in $X$.

Given a minimal parabolic subgroup $P_i$ and the closure $\bar\cO$
of an orbit, consider the morphism
$$
F_{i}\colon P_i\times{}^B \bar\cO\to P_i\bar\cO,
$$
$$
(p,x)\mapsto px.
$$
Suppose that $\bar \cO\neq P_i\bar\cO$. Knop \cite{K} and
Richardson--Springer \cite{RS} showed that the following three cases
can occur:

\begin{itemize}
\item Type U: $P_i\cO=\cO'\sqcup\cO$, and $F_i$ is
birational.
\item Type N: $P_i\cO=\cO'\sqcup\cO$, and $F_i$ is
of degree 2.
\item Type T: $P_i\cO=\cO'\sqcup\cO\sqcup\cO''$, and
$F_i$ is birational. In this case, $\dim\cO''=\dim\cO$.
\end{itemize}

It turns out that type case N never occurs in our situation.

\begin{prop}\label{onestepbirational} Let $\cO$ be a $B$-orbit in
$X$, and let $P_i$ be a minimal parabolic subgroup raising this
orbit. Then the map $F_i\colon P_i\times{}^B\cO\to P_i\cO$
 is birational.
\end{prop}

\proof Choose the canonical representative $x\in\cO$ as in
Prop.~\ref{baseschoice}. A straightforward calculation shows that
the isotropy group of $x$ in $P_i$ is equal to the stabilizer of $x$
in $B$, described in Prop.~\ref{isotr}. This implies the
birationality of $F_i$.\endproof

\remark The two remaining cases correspond to the two possible
``raisings" described in the subsection~\ref{combdescr}: (T)
corresponds to item 1, and (U) corresponds to item 2. In the first
case, the rank of the orbit is increased by one, and in the second
case, it does not change. Thus, the weak order is compatible with
the rank function: if $\cO\preceq\cO'$, then $\rk\cO\leq\rk\cO'$.
This is true in general for spherical varieties (cf., for instance,
\cite{B1}). Note that the strong order is \emph{not} compatible with
the rank function.

Proposition~\ref{onestepbirational}, together with
Theorem~\ref{minorbits}, allows us to construct desingularizations
for $\bar\cO$'s similar to Bott--Samelson desingularizations of
Schubert varieties in Grassmannians.

Given an orbit $\cO$, consider a minimal orbit $\cO_{min}$ that is
less than or equal to than $\cO$ with respect to the weak order.
This means that there exists a sequence $(P_{i_1},\dots,P_{i_r})$ of
minimal parabolic subgroups such that
$$
\bar\cO=P_{i_r}\dots P_{i_1}\bar\cO_{min}.
$$
Consider the map
$$
F\colon P_{i_r}\times^B \dots \times^B
P_{i_1}\times^B\bar\cO_{min}\to\bar\cO,
$$
$$
F\colon  (p_{i_r},\dots,p_{i_1},x)\mapsto p_{i_r}\dots p_{i_1}x.
$$
By Proposition~\ref{onestepbirational}, it is birational. But this
is not yet a desingularization, because $\bar\cO_{min}$ can be
singular.

The second step of desingularization consists in constructing a
$B$-equivariant desingularization for $\bar\cO_{min}$. We have
already proved in Theorem~\ref{minorbits} that $\bar\cO_{min}$ can
be represented as the direct product
$$
\bar\cO_{min}=X_w\times X_v
$$
for some Schubert varieties $X_w\subset\Gr(k,V)$ and
$X_v\subset\Gr(l,V)$.

For $X_w$ and $X_v$, one can take Bott--Samelson desingularizations
$$
F_w\colon Z_{w}\to X_w\quad \text{and}\quad F_v\colon Z_{v}\to X_v.
$$
(Details can be found, for instance, in \cite{B2}.) Thus, we obtain
the desingularization
$$
F_w\times F_v\colon\quad Z_{w}\times Z_{w}\to X_w\times
X_v=\bar\cO_{min}.
$$

Having this, we can combine this map with the map $F$ and obtain the
main result of this paper:
\begin{theorem} The map
$$
\tilde F=F\circ(F_w\times F_v)\colon\quad P_{i_r}\times^B \dots
\times^B P_{i_1}\times^B(Z_w\times Z_v)\to\bar\cO
$$
is a desingularization of $\bar\cO$.
\end{theorem}

\proof We have already seen that both maps $F$ and $F_w\times F_v$
are birational morphisms. Since all the considered varieties are
projective, these morphisms are proper. The variety $P_{i_r}\times^B
\dots \times^B P_{i_1}\times^B(Z_w\times Z_v)$ is a sequence of
homogeneous $B$-bundles with nonsingular fibers, and hence it is
nonsingular itself.\endproof

{\noindent\footnotesize{\sc Independent University of Moscow, Bolshoi Vlasievskii per., 11, 119002 Moscow, Russia\\
Institut Fourier, 100 rue des Maths, 38400 Saint-Martin d'H\`eres,
France}\\
\emph{E-mail address:} \verb"smirnoff@mccme.ru" }

\noindent{{\bf Keywords:} Grassmannian, spherical variety,
desingularization}

\end{document}